\documentclass[12pt]{article}
\usepackage{amsmath,amsfonts,amsthm,mathrsfs}
\newcommand{\R}{{\mathbb R}} \newcommand{\N}{{\mathbb N}}
\newcommand{\K}{{\mathbb K}} 
  
\newcommand{\Prm}{{\mathbb P}}

\renewcommand{\epsilon}{\varepsilon } 
 
\renewcommand{\rho}{\varrho } 
 
\renewcommand{\phi}{\varphi }

\newcommand{\E}{{\mathbb E}\,}
\newcommand{\ES}{{\mathbb E}'\,}

\newcommand{\ran}{{\rm ran }}
\newcommand{\de}{{\rm det }}
\newcommand{\ca}{{\rm card}}

\newtheorem{theorem}{Theorem}%[section]
\newtheorem{lemma}{Lemma}
\newtheorem{corollary}{Corollary}
\newtheorem{proposition}{Proposition}

\begin {document}
 \title{
 On the randomized complexity of Banach space valued integration}

\author  {Stefan Heinrich\\
Department of Computer Science\\
University of Kaiserslautern\\
D-67653 Kaiserslautern, Germany\\
e-mail: heinrich@informatik.uni-kl.de\\
\\
Aicke Hinrichs\\
Institute of Mathematics\\
University of Rostock\\
D-18051 Rostock, Germany\\
e-mail: aicke.hinrichs@uni-rostock.de
\\\\
%{\it Dedicated to Albrecht Pietsch on the occasion of his 80th birthday}
} 
\date{}
\maketitle
\begin{abstract}
We study the complexity of Banach space valued integration in the randomized setting. We are concerned with $r$-times continuously differentiable functions on the $d$-dimensional unit cube $Q$, with values in a Banach space $X$, and investigate the relation of the optimal convergence rate to the geometry of $X$. It turns out that the $n$-th minimal errors are bounded by $cn^{-r/d-1+1/p}$ if and only if $X$ is of equal norm type $p$.  
\end{abstract}
\section{Introduction}
\label{sec:1}
Integration of scalar valued functions is an intensively studied topic in the theory of information-based complexity, see \cite{TWW88}, \cite{Nov88}, \cite{NW10}. Motivated by applications to parametric integration, recently the complexity of Banach space valued integration was considered in \cite{DH12}.
It was shown that the behaviour of the $n$-th minimal 
errors $e_n^\ran$ of randomized integration in $C^r(Q,X)$ is related to the geometry of the Banach space $X$ in the following way: The infimum of the exponents of the rate is determined by the supremum of $p$ such that $X$ is of type $p$. In the present paper we further investigate this relation. We establish a connection between $n$-th minimal errors and equal norm type $p$ constants for $n$ vectors. It follows that $e_n^\ran$ is bounded by $cn^{-r/d-1+1/p}$ if and only if $X$ is of equal norm type $p$.  
\section{Preliminaries}
\label{sec:2}
Let ${\mathbb N}=\{1,2,\dots\}$ 
and  ${\mathbb N}_0=\{0,1,2,\dots\}$.
We introduce some notation and concepts from Banach space theory needed in the sequel. For Banach spaces $X$ and $Y$ let $B_X$ be the closed unit ball of $X$  and  
$\mathscr{L}(X,Y)$ the space of bounded linear operators from 
$X$ to $Y$, endowed with the usual norm.
If $X=Y$, we write $
\mathscr{L}(X)$. The norm of $X$ is denoted by $\| \cdot\|$, while other norms are distinguished by subscripts. We assume that all considered Banach spaces are defined over 
the same scalar field ${\mathbb K}={\mathbb R}$ or ${\mathbb K}={\mathbb C}$.

Let $Q=[0,1]^d$ and let $C^{r}(Q,X)$ 
be the space of all $r$-times continuously differentiable functions $f:Q\to X$ equipped with the norm
$$
\|f\|_{C^r(Q,X)}=\max_{0\le |\alpha|\le r,\, t\in Q}\| D^\alpha f(t)\|,
$$
where $\alpha=(\alpha_1,\dots,\alpha_d)$, $|\alpha|=|\alpha_1|+\dots+|\alpha_d|$ and  $D^\alpha$ denotes the respective partial derivative.
For $r=0$ we write $C^0(Q,X)=C(Q,X)$, which is the space of continuous $X$-valued functions on $Q$. If $X={\mathbb K}$, we write $C^r(Q)$ and $C(Q)$.

Let $1\le p\le 2$. A Banach space $X$ is said to be of (Rademacher) type $p$, if there is a constant $c>0$ such that for all $n\in{\mathbb N}$ and $x_1,\dots,x_n\in X$ 
\begin{equation}
\label{A6}
\left(\E\Big\| \sum_{i=1}^n \epsilon_i x_i\Big\|^p\right)^{1/p}\le c\left(\sum_{k=1}^n \| x_i\|^p\right)^{1/p},
\end{equation}
where $(\varepsilon_i)_{i=1}^n$ is a sequence of independent Bernoulli random variables with ${\mathbb P}\{\varepsilon_i=-1\}={\mathbb P}\{\varepsilon_i=+1\}=1/2$ on some probability space $(\Omega,\Sigma,\Prm)$ (we refer to \cite{MP76,LT91} for this notion and related facts). The smallest constant satisfying (\ref{A6}) is called the type $p$ constant of $X$ and is denoted by $\tau_p(X)$. If there is no such $c>0$, we put $\tau_p(X)=\infty$.
The space $L_{p_1}(\mathcal{N},\nu)$  with $(\mathcal{N},\nu)$  an arbitrary measure  space and $p_1<\infty$ is of type $p$ with $p=\min(p_1,2)$. 

Furthermore, given $n\in \N$, let $\sigma_{p,n}(X)$ be the smallest $c>0$ for which (\ref{A6}) holds for any $x_1,\dots,x_n\in X$ with $\|x_1\|=\dots=\|x_n\|$.
The contraction principle for Rademacher series, see (\cite{LT91}, Th.\ 4.4), implies that $\sigma_{p,n}(X)$ is the smallest constant $c>0$ such that for $x_1,\dots,x_n\in X$ 
\begin{equation}
\label{C3}
\left(\E\Big\| \sum_{i=1}^n \epsilon_i x_i\Big\|^p\right)^{1/p}\le cn^{1/p} \max_{1\le i\le n} \|x_i\|.
\end{equation}
We say that $X$ is of equal norm type $p$, if there is a constant $c>0$ such that $\sigma_{p,n}(X)\le c$ for all $n\in\N$. Clearly, 
$\sigma_{p,n}(X)\le \tau_p(X)$ and type $p$ implies equal norm type $p$.

Let us comment a little more on the relation of the different notions of type which are used here and in the literature.
The concept of equal norm type $p$ was first introduced and used by R. C. James in the case $p=2$ in \cite{Jam78}.
There it is shown that $X$ is of equal norm type 2 if and only if $X$ is of type 2.
This result is attributed to G. Pisier.
Later, it even turned out in \cite{BKT89} that the sequence $\sigma_{2,n}(X)$ and the corresponding
sequence $\tau_{2,n}(X)$ of type 2 constants computed with $n$ vectors are uniformly equivalent.
In contrast, for $1<p<2$, L. Tzafriri \cite{Tza79} constructed Tsirelson spaces without type $p$ but with equal norm type $p$. 
Finally, V. Mascioni introduced and studied the notion of weak type $p$ for $1<p<2$ in \cite{Mas88} and
showed that, again in contrast to the situation for $p=2$, a Banach space $X$ is of weak type $p$ if and only if
it is of equal norm type $p$.   

Throughout the paper $c,c_1,c_2,\dots$ are
constants, which depend only on the problem parameters $r,d$,
but depend neither on the algorithm parameters $n,l$ etc. nor on the input  $f$. The same symbol may denote different constants, 
even in a sequence of relations. 

For $r,k\in {\mathbb N}$ we let $P^{r,X}_k\in \mathscr{L}(C(Q,X))$ be $X$-valued composite tensor product Lagrange interpolation of degree $r$ with respect to the partition of $[0,1]^d$ into $k^d$ subcubes of sidelength $k^{-1}$ of disjoint interior, see \cite{DH12}.
 Given $r\in{\mathbb N}_0$ and $d\in {\mathbb N}$, there are constants $c_1,c_2>0$ such that 
for all Banach spaces $X$ and all $k\in{\mathbb N}$ 
\begin{eqnarray}
\sup_{f\in B_{C^r(Q,X)}}\|f-P^{r,X}_kf \|_{C(Q,X)}\le c_2k^{-r}\label{A3}
\end{eqnarray}
(see \cite{DH12}).
\section{Banach space valued  integration}
\label{sec:3}
Let $X$ be a Banach space, $r\in{\mathbb N}_0$, and let the integration operator $S^X:C(Q,X)\to X$ be given by 
\begin{equation*}
S^Xf=\int_Q f(t)dt.
\end{equation*}
We will work in the setting of information-based complexity theory, see \cite{TWW88,Nov88,NW10}. 
Below $e_n^{\rm det }(S^X,B_{C^r(Q,X)})$ and $e_n^\ran(S^X,B_{C^r(Q,X)})$ denote the $n$-th minimal error of $S^X$ on $B_{C^r(Q,X)}$ in the deterministic, respectively randomized setting, that is, the minimal possible error among all deterministic, respectively randomized algorithms, approximating $S^X$ on $B_{C^r(Q,X)}$ that use 
at most $n$ values of the input function $f$. The precise notions are recalled in the appendix.
The following was shown in \cite{DH12}.

\begin{theorem}
\label{theo:1}
Let $r\in \mathbb{N}_0$ and $1\le p\le 2$.  Then there are constants
$c_{1-4}>0$ such that for all Banach spaces $X$ and $n\in{\mathbb N}$ the following holds. The deterministic  $n$-th minimal error satisfies     
$$
c_1n^{-r/d}\le e_n^{\rm det }(S^X,B_{C^r(Q,X)})\le c_2n^{-r/d}.
$$
Moreover, if $X$ is of type $p$ and $p_X$ is the supremum of all $p_1$ such that $X$ is of type $p_1$,  
then  the randomized $n$-th minimal error   fulfills                                                       
$$
c_3n^{-r/d-1+1/p_X}\le e_n^{\rm ran }(S^X,B_{C^r(Q,X)})\le c_4\tau_p(X)n^{-r/d-1+1/p}.
$$
\end{theorem} 
As a consequence, we obtain 
\begin{corollary}
\label{cor:1}
Let $r\in {\mathbb N}_0$ and $1\le p\le 2$. Then the following are equivalent:\\[.1cm]
(i) $X$ is of type $p_1$ for all $p_1<p$.\\[.1cm]
(ii) For each $p_1<p$ there is a constant $c>0$ such that for all $n\in \N$ 
$$
e_n^{\rm ran }(S^X,B_{C^r(Q,X)})\le cn^{-r/d-1+1/p_1}.
$$
\end{corollary}
The main result of the present paper is the following
\begin{theorem} 
\label{theo:2}
Let $1\le p\le 2$ and $r\in\N_0$. Then there are constants $c_1,c_2>0$ such that  for all Banach spaces $X$ and all $n\in\N$
\begin{equation}
\label{A2}
c_1 n^{r/d+1-1/p}e_n^\ran(S^X,B_{C^r(Q,X)})\le \sigma_{p,n}(X)\le c_2\max_{1\le k\le n} k^{r/d+1-1/p}e_k^\ran(S^X,B_{C^r(Q,X)}).
\end{equation}
\end{theorem}
This allows to sharpen Corollary \ref{cor:1} in the following way.
\begin{corollary}
\label{cor:2}
Let $r\in {\mathbb N}_0$ and $1\le p\le 2$. Then the following are equivalent:\\[.1cm]
(i) $X$ is of equal norm type $p$.\\[.1cm]
(ii) There is a constant $c>0$ such that for all $n\in \N$ 
$$
e_n^{\rm ran }(S^X,B_{C^r(Q,X)})\le cn^{-r/d-1+1/p}.
$$
\end{corollary}
Recall from the preliminaries that the conditions in the corollary are also equivalent to\\[.1cm]
\emph{(iii) $X$ is of type $2$ if $p=2$ and of weak type $p$ if $1<p<2$, respectively.}\\[.1cm]

For the proof of Theorem \ref{theo:2} we need a number of auxiliary results.
The following lemma is a slight modification of Prop.\ 9.11 of \cite{LT91}, with essentially the same proof, which we include for the sake of completeness.
\begin{lemma}
\label{lem:2} 
Let $1\le p\le 2$. Then there is a constant $c>0$ such that for each Banach space  $X$, each $n\in \N$
 and each sequence of independent, essentially bounded, mean zero $X$-valued random variables  
$(\eta_i)_{i=1}^n$ on some probability space $(\Omega,\Sigma,\Prm)$ the following holds:
$$
\Bigg(\E \Big\|\sum_{i=1}^n \eta_i\Big\|^{p}\Bigg)^{1/p}
\le c\sigma_{p,n}(X)n^{1/p}\max_{1\le i\le n}
\|\eta_i\|_{L_\infty(\Omega,\Prm,X)}.
$$
\end{lemma}

\begin{proof}
Let $(\epsilon_i)_{i=1}^{n}$ be independent, symmetric Bernoulli random variables on some probability space
$(\Omega',\Sigma',\Prm')$ different from $(\Omega,\Sigma,\Prm)$. Considering $(\eta_i)_{i=1}^n$ and $(\epsilon_i)_{i=1}^{n}$ as random variables on the product probability space, we denote the expectation with respect to 
$\Prm'$ by $\ES$ (and the expectation with respect to $\Prm$, as before, by $\E$). Using Lemma 6.3 of \cite{LT91} and (\ref{C3}), we get
\begin{eqnarray*}
\Bigg(\E \Big\|\sum_{i=1}^n \eta_i\Big\|^{p}\Bigg)^{1/p}&\le& 
2\left(\E\ES\Big\|\sum_{i=1}^{n}\epsilon_i \eta_i\Big\|^{p}\right)^{1/p}
\\
&\le&2\sigma_{p,n}(X)n^{1/p}\left(\E\max_{1\le i\le n}\|\eta_i\|^p\right)^{1/p}
\nonumber\\
&\le&  2\sigma_{p,n}(X)n^{1/p}\max_{1\le i\le n}
\|\eta_i\|_{L_\infty(\Omega,\Prm,X)}.\quad
\label{AV2}
\end{eqnarray*}
\end{proof}

Next we introduce an algorithm for the aproximation of $S^Xf$. 
Let $n\in{\mathbb N}$ and let $\xi_i:\Omega\to Q$ $(i=1,\dots,n)$ 
be independent random variables 
on some  probability space $(\Omega,\Sigma,{\mathbb P})$
uniformly distributed on $Q$. 
Define for $f\in C(Q,X)$
\begin{equation}
\label{A8}
A^{0,X}_{n,\omega}f=\frac{1}{n}\sum_{i=1}^n f(\xi_i(\omega))
\end{equation}
and, if $r\ge 1$,  put $k=\left\lceil n^{1/d}\right\rceil$ and
\begin{equation}
\label{A9}
A^{r,X}_{n,\omega}f=S^X(P_k^{r,X}f)+A^{0,X}_{n,\omega}(f-P_k^{r,X}f).
\end{equation}
These are the Banach space valued versions of the standard Monte Carlo method ($r=0$) and the Monte Carlo method with separation of the main part ($r\ge 1$).
The following extends the second part of Proposition 1 of \cite{DH12}.
\begin{proposition}
\label{pro:2}
Let $r\in{\mathbb N}_0$ and $1\le p\le 2$. Then there is a constant $c>0$ such that for all Banach spaces $X$, $n\in{\mathbb N}$, and $f\in C^r(Q,X)$
\begin{eqnarray}
\left({\mathbb E}\, \|S^Xf-A_{n,\omega}^{r,X}f \|^p\right)^{1/p}
&\le& c \sigma_{p,n}(X)n^{-r/d-1+1/p}\|f\|_{C^r(Q,X)}.\label{C4}
\end{eqnarray}
\end{proposition}
\begin{proof}
Let us first consider the case $r=0$. Let $f\in C(Q,X)$ and put
$$
\eta_i(\omega)= \int_Q f(t)dt-f(\xi_i(\omega)).
$$
Clearly, ${\mathbb E}\,\eta_i(\omega) =0$,
$$
S^Xf-A_{n,\omega}^{0,X}f=\frac{1}{n}\sum_{i=1}^n\eta_i(\omega)
$$
and 
$$
\|\eta_i(\omega)\|\le 2\|f\|_{C(Q,X)}.
$$
An application of Lemma \ref{lem:2} gives (\ref{C4}).  
If $r\ge 1$,
we have
$$
 S^Xf-A^{r,X}_{n,\omega}f=S^X(f-P_k^{r,X}f)-A^{0,X}_{n,\omega}(f-P_k^{r,X}f) 
$$
and the result follows from (\ref{A3}) and the case $r=0$.
\end{proof}

\begin{lemma}
\label{lem:3} 
Let $1\le p\le 2$. Then there are constants $c>0$ and $0<\gamma<1$ such that for each Banach space  $X$, each $n\in \N$, and $(x_i)_{i=1}^n\subset X$
there is a subset $I\subseteq \{1,\dots,n\}$ with $|I|\ge \gamma n$ and
$$
\E \Big\|\sum_{i\in I} \epsilon_ix_i\Big\|
\le cn^{1/p}
\|(x_i)\|_{\ell_\infty^n(X)}\max_{1\le k\le n} k^{r/d+1-1/p}e_k^\ran(S^X,B_{C^r(Q,X)}).
$$
\end{lemma}
\begin{proof}
Since for $n\in \N$
$$
\max_{1\le k\le n} k^{r/d+1-1/p}e_k^\ran(S^X,B_{C^r(Q,X)})\ge e_1^\ran(S^\K,B_{C^r(Q,\K)})>0,
$$
the statement is trivial for $n<8^d$. Therefore we can assume $n\ge 8^d$. Clearly, we can also assume
$
\|(x_i)\|_{\ell_\infty^n(X)}>0.
$
Let $m\in {\mathbb N}$ be such that 
\begin{equation}
\label{A5}
m^d\le n< (m+1)^d,
\end{equation}
hence 
\begin{equation}
\label{B8}
m\ge 8.
\end{equation}
Let $\psi$ be an infinitely differentiable function on $\R^d$ such that $\psi(t) > 0$ for $t\in (0,1)^d$ and 
${\rm supp\,}\psi \subset [0,1]^d$. Let $(Q_i)_{i=1}^{m^d}$ be the partition of $Q$ into closed cubes of side length $m^{-1}$ of disjoint interior, let $t_i$ be the point in $Q_i$ with minimal coordinates and 
define $\psi_i\in C(Q)$ by
$$
\psi_i(t)=\psi(m(t-t_i))\quad (i=1, \dots,m^d).
$$
It is easily verified that there is a constant $c_0>0$  such that for all $(\alpha_i)_{i=1}^{m^d}\in [-1,1]^{m^d}$
$$
\Big\|\sum_{i=1}^{m^d}\alpha_i x_i\psi_i\Big\|_{C^r(Q,X)}\le c_0 m^r\|(x_i)\|_{\ell_\infty^n(X)}.
$$
Setting 
$$
f_i=c_0^{-1}m^{-r}\|(x_i)\|_{\ell_\infty^n(X)}^{-1}x_i\psi_i 
$$ 
it follows that
$$
\sum_{i=1}^{m^d}\alpha_i f_i\in B_{C^r(Q,X)} \qquad \mbox{for all}  \quad (\alpha_i)_{i=1}^{m^d}\in [-1,1]^{m^d}.
$$
Moreover, with $\sigma=\int_Q \psi(t) dt$ we have
\begin{eqnarray*}
&&\Big\|\sum_{i=1}^{m^d}\alpha_iS^Xf_i\Big\|
= c_0^{-1}m^{-r}\|(x_i)\|_{\ell_\infty^n(X)}^{-1}\Big\|\sum_{i=1}^{m^d}\alpha_ix_i\int_Q\psi_i(t)dt\Big\|\\
&=& c_0^{-1}\sigma m^{-r-d}\|(x_i)\|_{\ell_\infty^n(X)}^{-1}\Big\|\sum_{i=1}^{m^d}\alpha_ix_i\Big\|.
\end{eqnarray*}
Next we use Lemma 5 and 6 of \cite{Hei05a} with $K=X$ (although stated for $K={\mathbb R}$, Lemma 6 is easily seen to hold for $K=X$, as well) to obtain for all $l\in\N$ with $l< m^d/4$
\begin{eqnarray*}
e_l^{\rm ran }(S^X,B_{C^r(Q,X)})&\ge& \frac{1}{4}\min_{I\subseteq \{1,\dots, m^d\}, |I|\ge m^d-4l}{\mathbb E}\, \Big\|  \sum_{i\in I}\varepsilon_iS^Xf_i\Big\|\\
&\ge& cm^{-r-d}\|(x_i)\|_{\ell_\infty^n(X)}^{-1}{\mathbb E}\, \Big\|  \sum_{i\in I}\varepsilon_ix_i\Big\|.
\end{eqnarray*}
We put $l= \lfloor m^d/8\rfloor$. Then 
\begin{equation}
\label{B9}
m^d/16< l\le m^d/8.
\end{equation}
Indeed, by (\ref{B8}) the left-hand inequality clearly holds for $m^d<16$, while for $m^d\ge 16$ we get $ \lfloor m^d/8\rfloor>m^d/8-1\ge m^d/16$. We conclude that there is an $I\subseteq \{1,\dots, m^d\}$ with $|I|\ge m^d-4l\ge  m^d/2$ and
\begin{eqnarray*}
{\mathbb E}\, \Big\|  \sum_{i\in I}\varepsilon_ix_i\Big\|
&\le&
cm^{r+d}\|(x_i)\|_{\ell_\infty^n(X)}e_l^{\rm ran }(S^X,B_{C^r(Q,X)}) \\
&\le&cm^{r+d}l^{-r/d+1/p-1}\|(x_i)\|_{\ell_\infty^n(X)}\max_{1\le k\le n} k^{r/d+1-1/p}e_k^\ran(S^X,B_{C^r(Q,X)}) \\
&\le& cn^{1/p}
\|(x_i)\|_{\ell_\infty^n(X)}\max_{1\le k\le n} k^{r/d+1-1/p}e_k^\ran(S^X,B_{C^r(Q,X)}),
\end{eqnarray*}
where we used (\ref{A5}) and (\ref{B9}). 
Finally, (\ref{A5}) and (\ref{B8}) give
$$
|I|\ge  m^d/2\ge \frac{m^d}{2(m+1)^d}\,n\ge \frac{8^d}{2\cdot9^d}\,n.
$$

\end{proof}
\noindent {\it Proof of Theorem \ref{theo:2}.}
The left-hand inequality of (\ref{A2}) follows directly from Proposition \ref{pro:2}, since the number of function values involved in $A_{n,\omega}^{r,X}$ is bounded by $ck^d+n\le cn$, see also (\ref{D2}). 

To prove the right-hand inequality of (\ref{A2}), let $n\in \N$ and $x_1,\dots,x_n\in X$.
We construct by induction a partition of $K=\{1,\dots,n\}$ into a sequence of disjoint subsets $(I_l)_{l=1}^{l^*}$ 
such that for $1\le l\le l^*$
\begin{equation}
\label{B3}
|I_l|\ge \gamma \,\Big|K\setminus \bigcup_{j<l}I_j \Big| 
\end{equation}
and 
\begin{eqnarray}
\label{B4}
\lefteqn{\E\, \Big\|  \sum_{i\in I_l}\varepsilon_ix_i\Big\|}\nonumber\\
&\le& c\Big|K\setminus \bigcup_{j<l}I_j \Big|^{1/p}
\|(x_i)\|_{\ell_\infty^n(X)}\max_{1\le k\le n} k^{r/d+1-1/p}e_k^\ran(S^X,B_{C^r(Q,X)}),
\quad
\end{eqnarray}
where $c $ and $\gamma$ are the constants from Lemma \ref{lem:3}. For $l=1$ the existence of an $I_1$ satisfying (\ref{B3}--\ref{B4}) follows directly from Lemma \ref{lem:3}. Now assume that we already have a sequence of disjoint subsets $(I_l)_{l=1}^m$ of $K$ satisfying (\ref{B3}--\ref{B4}). If
$$
J:=K\setminus \bigcup_{j\le m}I_j\ne \emptyset,
$$
we apply Lemma \ref{lem:3} to 
$(x_i)_{i\in J}$ to find $I_{m+1}\subseteq J$ with 
\begin{equation}
\label{B5}
|I_{m+1}|\ge \gamma |J|
\end{equation}
and 
\begin{eqnarray}
\label{B6}
\lefteqn{\E\, \Big\|  \sum_{i\in I_{m+1}}\varepsilon_ix_i\Big\|}\nonumber\\
&\le& c|J|^{1/p}
\|(x_i)_{i\in J}\|_{\ell_\infty(J,X)}\max_{1\le k\le |J|} k^{r/d+1-1/p}e_k^\ran(S^X,B_{C^r(Q,X)}).
\end{eqnarray}
Observe that for $l=m+1$, (\ref{B5}) is just (\ref{B3}) and (\ref{B6}) implies (\ref{B4}).  
Furthermore, (\ref{B3}) implies
\begin{equation*}
\Big|K\setminus \bigcup_{j\le l}I_j \Big|\le (1-\gamma) \,\Big|K\setminus \bigcup_{j\le l-1}I_j \Big| 
\end{equation*}
and therefore 
\begin{equation}
\label{B7}
\Big|K\setminus \bigcup_{j\le l}I_j \Big|\le (1-\gamma)^l n. 
\end{equation}
It follows that the process stops with $K=\bigcup_{j\le l}I_j$ for a certain $l=l^*\in\N$. This completes the construction.

Using the equivalence of moments (Theorem 4.7 of \cite{LT91}), we get from (\ref{B4}) and (\ref{B7})
\begin{eqnarray*}
\lefteqn{ \left(\E\Big\| \sum_{i=1}^n \epsilon_i x_i\Big\|^p\right)^{1/p}  }\\
&\le & c\,\E\Big\| \sum_{i=1}^n \epsilon_i x_i\Big\|\le c\sum_{l=1}^{l^*}\E\Big\| \sum_{i\in I_l} \epsilon_i x_i\Big\|\\
&\le &cn^{1/p}\|(x_i)\|_{\ell_\infty^n(X)}\max_{1\le k\le n} k^{r/d+1-1/p}e_k^\ran(S^X,B_{C^r(Q,X)})\sum_{l=1}^{l^*}(1-\gamma)^{(l-1)/p}.
\end{eqnarray*}
This gives the upper bound of (\ref{A2}). 

\qed
\medskip

Let us mention that results analogous to Theorem \ref{theo:2} and Corollary \ref{cor:2} above also hold for Banach space valued indefinite integration (see \cite{DH12} for the definition) and for the solution of initial value problems for Banach space valued ordinary differential equations \cite{Hei12}. Indeed, an inspection of the respective proofs together with Lemma \ref{lem:2} of the present paper shows that Proposition 2 of \cite{DH12} also holds with $\tau_p(X)$ replaced by $\sigma_{p,n}(X)$, and similarly Proposition 3.4 of \cite{Hei12}. Moreover, in both papers the lower bounds on $e_n^\ran$ are obtained by reduction to (definite) integration and thus the righ-hand side inequality of (\ref{A2}) carries over directly.

\section{Appendix}
In this appendix we recall some basic notions  of 
information-based complexity -- the framework we used above. We refer to \cite{Nov88,TWW88} for more on this subject and to \cite{Hei05a, Hei05b} for the particular notation applied here.
First we introduce the class of deterministic adaptive algorithms of varying cardinality 
$
\mathcal{A}^{\det}(C(Q,X),X)
$.
It consists of tuples 
$
A=((L_i)_{i=1}^\infty, (\rho_i)_{i=0}^\infty,(\phi_i)_{i=0}^\infty),
$
with
$
L_1\in Q$, $\rho_0\in\{0,1\}$, $\phi_0\in X$
and 
$$
L_i : X^{i-1}\to Q\quad (i=2,3,\dots),\quad
 \rho_i:  X^i\to \{0,1\},\;
\phi_i:  X^i\to X \quad (i=1,2,\dots)
$$
being arbitrary mappings.
To each $f\in C(Q,X)$, we associate a sequence $(t_i)_{i=1}^\infty$ with $t_i\in Q$ as follows:
\begin{equation*}
t_1=L_1, \quad 
t_i=L_i(f(t_1),\dots,f(t_{i-1}))\quad(i\ge 2).
\end{equation*}
Define $\ca(A,f)$, the cardinality  of $A$ at input $f$, to be $0$ if $\rho_0=1$. If $\rho_0=0$, let $\ca(A,f)$ be
the first integer $n\ge 1$ with 
$
\rho_n(f(t_1),\dots,f(t_n))=1,
$
if there is such an $n$, and $\ca(A,f)=+\infty$ otherwise.
For $f\in C(Q,X)$ with $\ca(A,f)<\infty$ we define the output $Af$ of algorithm $A$ at input $f$ as
$$
Af=\left\{\begin{array}{lll}
 \phi_0  & \mbox{if} \quad n=0   \\
 \phi_n(f(t_1),\dots,f(t_n))  &\mbox{if} \quad n\ge 1.    
    \end{array}
\right.
$$
Let $r\in \N_0$. Given $n\in\N_0$, we let $\mathcal{A}_n^{\det}(B_{C^r(Q,X)},X)$ be the set of those
$A\in \mathcal{A}^{\det}(C(Q,X),X)$ for which
$$
\max_{f\in B_{C^r(Q,X)}} \,\ca (A,f)\le n.
$$
The error of $A\in \mathcal{A}_n^{\det}(B_{C^r(Q,X)},X)$
as an approximation of  $S^X$ is defined as
$$
e(S^X,A,B_{C^r(Q,X)})=\sup_{f\in B_{C^r(Q,X)}} \|S^Xf-Af\|.
$$
The deterministic $n$-th minimal error of $S^X$ is defined for $n\in\N_0$ as
\begin{equation*}
e_n^\de (S^X,B_{C^r(Q,X)})=\inf_{A\in\mathcal{A}_n^\de(B_{C^r(Q,X)}) }  
e(S^X,A,B_{C^r(Q,X)}).
\end{equation*}
It follows that no deterministic algorithm that uses at most $n$ function values can have  a smaller error 
than $e_n^\de (S^X,B_{C^r(Q,X)})$.

Next we introduce the class of randomized adaptive algorithms of varying cardinality
$
\mathcal{A}_n^{\ran}(B_{C^r(Q,X)},X),
$
consisting of tuples
$
A=((\Omega,\Sigma,\Prm),(A_\omega)_{\omega\in \Omega}),
$
where $(\Omega,\Sigma,\Prm)$ is a probability space,
$
A_\omega\in  \mathcal{A}^{\det}(C(Q,X),X)
$
for all $\omega\in \Omega$, 
and for each $f\in B_{C^r(Q,X)}$
the mapping 
$
\omega\in\Omega\to \ca (A_\omega, f)
$
is $\Sigma$-measurable and  satisfies
$
\E\,\ca (A_\omega, f)\le n.
$
Moreover, the mapping
$
\omega\in\Omega\to A_\omega f\in X 
$
is  $\Sigma$-to-Borel measurable and essentially separably valued,
 i.e., there is a separable subspace
$X_0\subseteq X$  such that 
$A_\omega f\in X_0$  for $\Prm$-almost all $\omega \in \Omega$.
The error of  $A\in\mathcal{A}_n^{\ran}(C(Q,X),X)$ in approximating
 $S^X$  on $B_{C^r(Q,X)}$ is defined as
\begin{equation*} 
e(S^X,A,B_{C^r(Q,X)}) = 
\sup_{f \in B_{C^r(Q,X)}}  \,\E  \| S^Xf- A_\omega f\|,
\end{equation*} 
and the randomized $n$-th minimal error of $S^X$ as
\begin{equation*}
e_n^\ran (S^X,B_{C^r(Q,X)})=\inf_{A\in\mathcal{A}_n^\ran(B_{C^r(Q,X)}) }  
e(S^X,A,B_{C^r(Q,X)}).
\end{equation*}
Consequently,  no randomized  
algorithm that uses (on the average) at most $n$ function values 
has an  error smaller
than $e_n^\ran (S^X,B_{C^r(Q,X)},X)$.

Define for $\epsilon>0$  the information complexity  as 
$$
n_\epsilon^\ran(S,B_{C^r(Q,X)})
=\min\{n\in\N_0:\,e_n^\ran(S,B_{C^r(Q,X)})\le \epsilon \},
$$
if there is such an $n$, and 
$
n_\epsilon^\ran(S,B_{C^r(Q,X)})=+\infty,
$
if there is no such $n$.
Thus, if $n_\epsilon^\ran(S,B_{C^r(Q,X)})<\infty$, it follows that any algorithm with error $\le \epsilon$ needs at least 
$n_\epsilon^\ran(S,B_{C^r(Q,X)})$ function values, while  $n_\epsilon^\ran(S,B_{C^r(Q,X)})=+\infty$
means that no algorithm at all has error 
$\le \epsilon$.
The information complexity is essentially the inverse function of the $n$-th minimal error. So determining the latter means determining the information complexity of the problem.

Let us also mention the subclasses consisting of quadrature formulas. Let $n\ge 1$. A mapping  $A:C(Q,X)\to X$ is called a deterministic quadrature formula with $n$ nodes, if there are $t_i\in Q$ and  $a_i\in \K$ ($1\le i\le n$) such that
$$
A f= \sum_{i=1}^n   \, a_if(t_i) \quad(f\in C(Q,X)). 
$$
In terms of the definition of  $\mathcal{A}^\de(C(Q,X),X)$ this means that the respective functions $L_i$ and $\rho_i$ are constant,  $
\rho_0=\rho_1=\dots=\rho_{n-1}=0$, $\rho_n=1$, and $\phi_n$ has the form $\phi_n(x_1,\dots,x_n)=\sum_{i=1}^n   a_ix_i$. Clearly, $A\in \mathcal{A}_n^{\de}(B_{C^r(Q,X)}, X)$. 

A tupel $A=((\Omega,\Sigma,\Prm),(A_\omega)_{\omega\in \Omega})$ is called a randomized quadrature with $n$ nodes if  there exist random variables $t_i:\Omega\to Q$ and  $a_i:\Omega\to \K$ ($1\le i\le n$) with 
$$
A_\omega f= \sum_{i=1}^n   \, a_i(\omega)f(t_i(\omega)) \quad(f\in C(Q,X),\,\omega\in \Omega). 
$$
For each such $A$ we have $A\in \mathcal{A}_n^{\ran}(B_{C^r(Q,X)}, X)$. Finally we note that 
the algorithms $A^{r,X}_{n,\omega}$ defined in (\ref{A8}) and (\ref{A9}) are 
quadratures. Indeed, for $A^{0,X}_{n,\omega}$ given by (\ref{A8}) this is obvious. For $r\ge 1$ we represent $P_k^{r,X}\in \mathscr{L}(C(Q,X))$ as
$$
P_k^{r,X}f= \sum_{j=1}^{M}f(u_j)\psi_j(t)
$$
with  $M\le ck^d$, $u_j\in Q$, $\psi_j\in C(Q)$ ($1\le i\le M$), and obtain, setting $b_j=\int_Q \psi_j(t)dt$,
\begin{eqnarray}
A^{r,X}_{n,\omega}f
&=&
S^X(P_k^{r,X}f)+A^{0,X}_{n,\omega}(f-P_k^{r,X}f)\nonumber\\
&=&
\sum_{j=1}^{M}b_jf(u_j)+\frac{1}{n}\sum_{i=1}^n\left(f(\xi_i(\omega))-\left(P_k^{r,X}f\right)(\xi_i(\omega))\right)\nonumber\\
&=&\sum_{j=1}^{M}b_jf(u_j)+\frac{1}{n}\sum_{i=1}^n f(\xi_i(\omega))-\sum_{j=1}^{M}\left(\frac{1}{n}\sum_{i=1}^n\psi_j(\xi_i(\omega))\right)f(u_j).\label{D2}
\end{eqnarray}

\end{document}